%% file: BlueCollianderWeightedNLSData.tex
\documentclass{amsart}
\usepackage{amssymb}
\usepackage{amsmath}
\usepackage{amsthm}

\usepackage{epsfig}
\usepackage{color}

\title[ Global existence for NLS weighted data.]{Global well-posedness in Sobolev space implies global existence for weighted $L^2$ initial data for $L^2$-critical NLS.}
\author{P. Blue, J. Colliander}
\date{\today}

\newcommand{\breakcompile}{\BreakCompile}

\newcommand{\JimsMacDoesntLikeFig}[1]{#1}

\newcommand{\oldversion}[1]{\breakcompile}

\newcommand{\preliminarycomments}[1]{\breakcompile}

\newcommand{\hide}[1]{\breakcompile}

\newcommand{\nls}{nonlinear Schr\"odinger }
\newcommand{\ls}{linear Schr\"odinger }

\newtheorem{theorem}{Theorem}[section]
\newtheorem{corollary}[theorem]{Corollary}
\newtheorem{lemma}[theorem]{Lemma}
\newtheorem{proposition}[theorem]{Proposition}
\newtheorem{definition}[theorem]{Definition}
\newtheorem{remark}[theorem]{Remark}

\numberwithin{equation}{section}
\numberwithin{theorem}{section}

\newtheorem{faketheoremaux}{Theorem}

\newtheorem{fakecorollaryaux}{Corollary}

\newcommand{\Reals}{\mathbb{R}}
\newcommand{\Complex}{\mathbb{C}}

\newcommand{\Fourier}[1]{\mathcal{F}[#1]}

\newcommand{\transT}{\tau}
\newcommand{\RegExponent}{s}
\newcommand{\Xsigma}{H^{0,\RegExponent}}
\newcommand{\Xone}{H^{0,1}}

\newcommand{\Jax}{\langle x\rangle} 

\newcommand{\Conf}{\mathcal{C}}

\newcommand{\Dt}{\frac{d}{dt}}
\newcommand{\dt}{\partial_t}

\newcommand{\deltacomesfromphrase}{coming from the local well-posedness theory }

\newcommand{\SBalancedExponent}{\frac{2(d+2)}{d}}

\begin{document}

\maketitle
\begin{abstract}
The $L^2$-critical defocusing nonlinear Schr\"odinger initial value problem on $\Reals^d$ is known to be locally well-posed for initial data in $L^2$. Hamiltonian conservation and the pseudoconformal transformation show that global well-posedness holds for initial data $u_0$ in Sobolev $H^1$ and for data in the weighted space $(1+|x|) u_0 \in L^2$. For the $d=2$ problem, it is known that global existence holds for data in $H^s$ and also for data in the weighted space $(1+|x|)^{\sigma} u_0 \in L^2$ for certain $s, \sigma < 1$. We prove: If global well-posedness holds in $H^s$ then global existence and scattering holds for initial data in the weighted space with $\sigma = s$. 
\end{abstract}

\section{Introduction}

Consider the initial value problem for the $L^2$ critical \nls equation for 
$u:\Reals\times\Reals^d\rightarrow\Complex$,
\begin{equation}
\label{NLS}
\left\{
\begin{matrix}
i\dt u+\Delta u=\lambda |u|^{\frac{4}{d}}u \\
u(t_0,x)= u_0(x) .
\end{matrix}
\right.
\end{equation}
This problem is called {\it defocusing} for $\lambda>0$ and {\it focusing} for $\lambda<0$. In dimension $d=2$, equation \eqref{NLS} reduces to the cubic \nls equation, $i\dt u+\Delta u=\lambda|u|^2u$, which appears widely as a model equation in Physics \cite{SulemSulem}. 

This problem is locally well-posed in $L^2$ or in any $H^\RegExponent$
with $\RegExponent\geq0$. That is, given initial data $u_0\in
H^\RegExponent(\Reals^d)$ with $\RegExponent\geq0$ there is a local existence
time, $T_{lwp}$, and a local in time solution $u:[t_0-T_{lwp},t_0 +T_{lwp}]\times\Reals^d\rightarrow\Complex$ such that $u$ solves \eqref{NLS} and the function $u \in C_t H^\RegExponent.$ For $\RegExponent>0$, $T_{lwp}$ is a decreasing function of the $H^\RegExponent$ norm of $u_0$. 

An open problem is to prove global well-posedness in $L^2$ in the defocusing case and under an appropriate smallness condition in the focusing case.  For the focusing case, it is believed that solutions with $L^2$ norm smaller than the ground state\footnote{The {\it{ground state}} is the unique (up to translations) positive solution of $-Q + \Delta Q = Q^3$.} mass $\| Q \|_{L^2}$ do not blowup and in fact scatter. Explicit blow-up solutions with Schwartz class initial data and with the mass of the ground state are known to exist in the focusing case. By finding the optimal constant in the Galiardo-Nirenberg estimate, $\|u\|_{L^\SBalancedExponent}^\SBalancedExponent \leq ( \frac{d+2}{d} \| Q \|^{-\frac4d}_{L^2}) \|\nabla u\|_{L^2}^{2} \|u\|_{L^2}^{\frac4d}$, Weinstein \cite{Weinstein} proved that $H^1$ initial data with $L^2$ norm less than the ground state mass evolves globally in time. In the defocusing case, $L^2$ solutions are expected to exist globally in time and scatter. Although the $L^2$ norm of $u(t)$ is constant on the local well-posedness time interval, this norm does not control the length of the local well-posedness time $T_{lwp}$ and can not be used to prove global well-posedness in $L^2$. In $H^1$, this problem has an additional conserved quantity, the energy, $E[u]=\int\frac12 |\nabla u|^2+\frac\lambda4|u|^4 dx$. In the defocusing case, the energy is positive and dominates the $H^1$ norm. Since the local well-posedness time is a function of the $H^1$ norm, at any time the solution persists for a uniformly long local well-posedness time, and, hence, globally in time. 

Bourgain was the first to prove \cite{Bourgain}, for the cubic problem in $\Reals^2$, global well-posedness below the energy
threshold $H^1$ by proving global well-posedness for data in
$H^\RegExponent$ in the defocusing case for $\RegExponent>\frac35$. The method in \cite{Bourgain} involved a decomposition of the
data into high and low frequencies with a sharp cut-off function in
the Fourier variables. Later, the ``$I$-method'' \cite{CKSTT} was used
to improve this to $\RegExponent>\frac47$ for the cubic problem on $\Reals^2$. 

The \nls equation \eqref{NLS} has a discrete pseudoconformal symmetry, 
\begin{equation}
\label{2.1}
\Conf[u](\transT,y)= v(\transT,y)=\frac{1}{|\transT|^\frac{d}{2}} e^\frac{iy^2}{4\transT} u(-\frac1\transT,\frac{y}{\transT}).
\end{equation}
This is a symmetry in the sense that, if $u(t,x)$ is a solution to the
\nls equation on $(t,x)\in[t_1,t_2]\times\Reals^d$, then
$v(\transT,y)$ is a solution on $\transT\in[-t_1^{-1},-t_2^{-1}], ~y
\in \Reals^d$. Throughout this paper, we shall use $u$, $t$, and $x$ to refer to a solution, the time variable, and the spatial variable respectively, to use $v$ as the pseudoconformal transform of $u$, and to use $\transT$ and $y$ as the arguments of $v$ which will be called the transformed time and space variables. In particular, $\transT=t^{-1}$. 

We introduce the space
\begin{equation}
\Xsigma=\{u\in L^2: |x|^\RegExponent u\in L^2\}
\end{equation}
with norm
\begin{equation}
\|u\|_{\Xsigma} = \|\Jax^\RegExponent u\|_{L^2}
\end{equation}
where $\Jax^2=1+x^2$. 

The pseudoconformal transform satisfies
\begin{equation}
\label{3.1}
E[v](\transT)=\| xu_0 \|_{L^2}^2 =\|u_0\|_{\Xone}^2 .
\end{equation}
The energy on the left side of this equation is already known to be
independent of $\transT$. This property was used to prove
global well-posedness in $L^2$ for initial data $u_0\in\Xsigma$
with $\RegExponent>\frac35$ in the defocusing case \cite{Bourgain}. The
proof involves a spatial decomposition analogous to the 
Fourier decomposition used in proving the $H^\RegExponent$ global
existence result. As a further consequence, \cite{Bourgain} also proves scattering, that is the existence of functions $u_\pm\in L^2$ such that
\begin{equation*}
\lim_{t\rightarrow\pm\infty} \| u(t) - e^{\pm i t\Delta} u_\pm\|_{L^2} =0.
\end{equation*}

This paper establishes that global well-posedness in $H^s$ for
\eqref{NLS} implies global existence and scattering in $L^2$ for
initial data in $H^{0,s}$ for \eqref{NLS}. Thus, the link between
$H^s$ global well-posedness and the evolution properties of $H^{0,s}$
initial data found in the $\Reals^2$ case in \cite{Bourgain} is in fact
common to all pseudoconformal or $L^2$-critical \nls initial value problems \eqref{NLS}.

\begin{proposition}
\label{P3.1}
Assume that the \nls equation \eqref{NLS} is globally well-posed in $H^\RegExponent$ (with the additional hypothesis that the initial data has $L^2$ norm bounded by $\| Q \|_{L^2}$ in the focusing case). 

If $u_0\in\Xsigma$ (and $\|u_0\|_{L^2}<\| Q \|_{L^2}$ in the focusing case), then there is a
function $u:\Reals\times\Reals^d\rightarrow\Complex$ which solves the
\nls equation \eqref{NLS} for all time. Furthermore,  there are functions $u_\pm\in \Xsigma$ such that 
\begin{equation*}
\lim_{t\rightarrow\pm\infty} \|e^{\mp it\Delta}u(t)-u_\pm\|_{\Xsigma} =0.
\end{equation*}
\end{proposition}

Based on the $H^\RegExponent$ global well-posedness result in \cite{CKSTT}, we obtain as a consequence of Proposition \ref{P3.1} that for $\RegExponent > \frac{4}{7}$, initial data in $\Xsigma$ evolve globally in time and scatter in $L^2$ under the cubic \nls flow on $\Reals^2$.

Global existence in $L^2$ for initial data in $\Xsigma$ means initial data $u_0\in\Xsigma$, which is also in $L^2$, evolves as a solution in $L^2$ and that this solution exists for all time. This is significantly weaker than global well-posedness in $H^\RegExponent$, which means that initial data in $H^\RegExponent$ continuously evolves in $H^\RegExponent$, that this solution exists for all time, and that the time evolution map $S_{NLS}(t,0):u_0\mapsto u(t)$ is continuous from $H^\RegExponent$ to $H^\RegExponent$. 

The asymmetry between $H^\RegExponent$ and $\Xsigma$ in Proposition \ref{P3.1} is a consequence of the local theory. The initial value problem \eqref{NLS} is locally well-posed in $H^\RegExponent$; whereas, in $\Xsigma$, \eqref{NLS} is ill-posed, so global existence in $L^2$ for initial data in $\Xsigma$ can not be extended to $\Xsigma$ global well-posedness. 

In the remainder of the introduction, we review $L^2$ and $H^\RegExponent$ local-well posedness, $\Xsigma$  ill-posedness, and some properties of the pseudoconformal transform. In Section \ref{s2}, we prove Proposition \ref{P3.1} by showing that, for a solution with initial data in $\Xsigma$, the pseudoconformal transform is in $H^\RegExponent$. This is done by taking regularized approximators and showing their transforms converge in $H^\RegExponent$ at a particular transformed time $-T_{lwp}^{-1}$. In Section \ref{s3}, we show that scattering is a consequence of the construction in Section \ref{s2}. 

We use the notation 
$S_{NLS}(t_2,t_1)$ to denote the \nls evolution map from time $t_1$ to time $t_2$, 
$\Fourier{\bullet}$ for the Fourier transform, 
and $\Re \alpha$ and $\Im \alpha$ for the real and imaginary part of $\alpha$ respectively.

\subsection{Local well-posedness theory}
\label{ssLWP}

The local well-posedness theory (see \cite{CazText},
\cite{LinaresPonce} for a review) begins with the presentation of the \nls equation as an integral equation through Duhamel's principle
\begin{align}
\label{5.1}
u=&\Phi_{u_0}[u] ,\\
\label{5.2}
\Phi_{u_0}[u]=&e^{it\Delta}u_0 + \int_0^t e^{i(t-t')\Delta}(\lambda|u|^2u)(t')dt'.
\end{align}
To prove that \eqref{5.1} has a unique solution, it is sufficient to show that $\Phi_{u_0}$ is a contraction in an appropriate space. This space will be the Strichartz space defined below. 

\begin{definition}
The pair $(q,r)$ is $L^2$ Strichartz admissible, or simply admissible, if 
\begin{align*}
\frac2q+\frac{d}{r}=& \frac{d}{2} ,
& 2\leq q\leq& \infty ,
& 2\leq r\leq& 2+\frac{4}{d-2} .
\end{align*}
In dimension $d=2$, there is the additional restriction that $2<q\leq\infty$ and $2\leq r<\infty$. In dimension $d=1$, $2\leq q\leq \infty$ and $2\leq r\leq\infty$. 

For $s \geq 0$, a pair $(q,r)$ is $H^\RegExponent$ Strichartz admissible if
\begin{align*}
\frac2q+\frac{d}{r}=&\frac{d}{2}-\RegExponent ,
& 2\leq q\leq& \infty , 
& 2\leq & r .
\end{align*}

The $L^2$ Strichartz norm, or simply the Strichartz norm, is
\begin{align*}
\|u\|_{S^0}=\sup_{\text{$(q,r)$ admissible}} \|u\|_{L^q_t L^r_x} .
\end{align*}

The homogeneous $H^\RegExponent$ Strichartz norm is 
\begin{align*}
\|u\|_{\dot{S}^\RegExponent}=\sup_{\text{$(q,r)$ admissible}} \|D^\RegExponent u\|_{L^q_t L^r_x} ,
\end{align*}
where $D$ is the Fourier multiplier defined by $\Fourier{ D f }(\xi) =
|\xi| \Fourier{f} (\xi)$. 

For $s > 0$, the $H^\RegExponent$ Strichartz norm is 
\begin{align*}
\|u\|_{S^\RegExponent}=\sup_{0 \leq \sigma \leq \RegExponent}
\|u\|_{\dot{S}^\sigma}
\end{align*}

For an interval $I$, the spaces $S^0(I)$, $\dot{S}^\RegExponent(I)$, and $S^\RegExponent(I)$ are the spaces with the above Strichartz norm where the $t$ integration is taken over the interval $t\in I$. 
\end{definition}

With this notation, we record the {\it{Strichartz estimates}}: For $s
> 0$, 
\begin{align*}
  \| e^{it \Delta} u_0 \|_{S^0} & \leq C \|u_0 \|_{L^2}, \\
\| e^{it \Delta} u_0 \|_{\dot{S}^s} & \leq C \| u_0 \|_{\dot{H}^s}, \\
\| e^{it \Delta} u_0 \|_{S^s } & \leq C \| u_0 \|_{H^s}.
\end{align*}

Using Duhamel's principle, estimates for solutions to the \ls equation, and a contraction argument, it is possible to prove local well-posedness, uniqueness, and continuity of solutions in the $L^2$ Strichartz space. The balanced pair $q=r=\SBalancedExponent$ plays an important role in our presentation of the local well-posedness theory. 

\begin{theorem}[$L^2$ local well-posedness]
\label{T6.1}
For $u_0\in L^2$ and $\delta_1>0$, there is a $T_{lwp}=T_{lwp}(u_0,\delta_1)>0$ and $u:[0,T_{lwp}]\times\Reals^d\rightarrow\Complex$ such that: 
\begin{enumerate}
\item $u$ solves the \nls equation \eqref{NLS} with initial data $u_0$ $\forall~ t\in[0,T_{lwp}]$;  
\item $\| u \|_{S^0_{([0,T_{lwp}])}} \leq 2 \|u_0 \|_{L^2}$; 
\item $\| u \|_{L^\SBalancedExponent_{tx}([0,T_{lwp}])} \leq \delta_1$;  
\item for all $t\in[0,T_{lwp}]$, $\|u(t)\|_{L^2}=\|u_0\|_{L^2}$; 
\item there is a $\delta_2>0$, such that, if $u_0'\in L^2$ with $\|u_0-u_0'\|_{L^2}<\delta_2$,
  then there is a $u'\in S^0([0,T_{lwp}])$ such that $u'$ solves the
  \nls equation with initial data $u_0'$; and $\|u-u'\|_{S^0([0,T_{lwp}])}\leq 2\|u_0-u_0'\|_{L^2}$. In particular, for $t\in[0,T_{lwp}]$, 
\begin{equation}
\|u(t)-u'(t)\|_{L^2}\leq 2\|u_0-u_0'\|_{L^2} .
\end{equation}
\end{enumerate}
\end{theorem}

The $T_{lwp}$ from this theorem will be referred to as the $L^2$ local well-posedness time. The $L^2$ maximal forward time of existence, $T^*$, is the time for which there is a solution $u:[0,T^*)\rightarrow L^2$ but no solution $\tilde{u}:[0,t]\rightarrow L^2$ for $t>T^*$.

The $L^2$ local well-posedness theory and
$L^\SBalancedExponent_{tx}([0,T_{lwp}])$ norm play a central role in
our presentation of the $H^\RegExponent$ theory. The $L^2$ local well-posedness time,
$T_{lwp}$, is the same as the $H^\RegExponent$ local well-posedness
time. If a solution has small $L^\SBalancedExponent_{tx}([0,T_{lwp}])$
norm, then its $H^\RegExponent$ norm can not grow by more than a small
factor, and, if another solution has initial data close in $L^2$, then
the $H^\RegExponent$ distance between the two solutions can not grow
by more than a small factor. The $L^2$ local well-posedness Theorem \ref{T6.1} uses $T_{lwp}$ chosen to enforce the $L^\SBalancedExponent_{tx}([0,T_{lwp}])$ smallness condition. 

\begin{theorem}[$H^\RegExponent$ local well-posedness]
\label{T7.1}
For $u_0\in H^\RegExponent$, the solution $u$ and local well-posedness time $T_{lwp}$ given in Theorem \ref{T6.1} satisfy: 
\begin{enumerate}
\item \label{HSigmaExistence} for $t\in[0,T_{lwp}]$, $u(t)\in H^\RegExponent$; 
\item \label{DeltaCHSigmaDoubling} there is a $\delta_3>0$, such that, if $\|u\|_{L^\SBalancedExponent_{tx}(I)}\leq \delta_3$, then $\|u\|_{\dot{S}^\RegExponent(I)}\leq 2 \|u_0\|_{H^\RegExponent}$; and 
\item \label{HSigmaDivergence}there is a $\delta_4>0$, such that, if $\|u\|_{L^\SBalancedExponent_{tx}(I)}\leq \delta_3$, $u'$ is a solution to the \nls equation with initial data $u'_0\in H^\RegExponent$, and $\|u_0-u'_0\|_{L^2}\leq\delta_4$, then $\|u-u'\|_{{S}^\RegExponent}\leq 2\|u_0-u'_0\|_{H^\RegExponent}$. In particular, 
\begin{equation}
\|u(t)-u'(t)\|_{H^\RegExponent} \leq 2 \|u_0-u'_0\|_{H^\RegExponent} .
\end{equation}
\end{enumerate}
\end{theorem}

Since $u$ can be extended to any interval on which the $L^2$ Strichartz norm is finite, the $L^2$ Strichartz norm must diverge on intervals approaching the maximal forward time of existence. 

\begin{corollary}[Maximal time blow up theorem]
\label{MaxTimeBlowUp}
If $u$ is a solution to the \nls equation \eqref{NLS} with initial data $u_0\in L^2$ and maximal forward time of existence $T^*$ (that is $u:[0,T^*)\rightarrow\Complex$ is a solution, but for $t>T^*$, there is no solution $\tilde{u}:[0,t]\rightarrow\Complex$), then 
\begin{equation}
\label{t*maximal}
\lim_{t\uparrow {T^*}} \| u \|_{S^0([0,t])}=\infty .
\end{equation}
\end{corollary}

\subsection{Local evolution of $\Xsigma$ initial data} 
\label{ssSpatiallyLocalised}

As noted earlier, Proposition \ref{P3.1} only asserts global existence in $L^2$ for initial data in $\Xsigma$, not $\Xsigma$ global well-posedness, because the \nls equation is ill-posed in $\Xsigma$. This ill-posedness is a consequence of the ill-posedness of the \ls equation in $\Xsigma$. It has solutions which start with finite initial $\Xsigma$ norm, but have divergent $\Xsigma$ norm at all later times. From these, it is possible to show that the \nls time evolution $S_{NLS}(t,0):u_0\mapsto u(t)$ is not continuous from $\Xsigma$ to $\Xsigma$. This result does not have an analogue in $H^\RegExponent, ~s>0$. 

Solutions to the \ls equation which have finite initial $\Xsigma$ norm but divergent $\Xsigma$ norm can be constructed from sums of Gaussians. Since the Fourier variable corresponds to velocity for the \ls equation, Gaussian initial data of width $a^{-\frac12}$ has a Gaussian Fourier transform with a width $a^\frac12$, and hences disperses linearly in time at speed $a^\frac12$. A linear combination of such initial data can be constructed to disperse with unbounded velocity and instantaneously divergent weighted norm. 

\begin{lemma}
\label{LM3.1}
For $s>0$, there is a solution $\Psi$ to the \ls equation with initial data $\Psi_0\in\Xsigma$ for which $\Psi(t)$ is not in $\Xsigma$ at any future time $t>0$. 
\end{lemma}

\begin{proof}
Let
\begin{align*}
(C_{d,s})^2=&\frac{1}{\pi^\frac{d}{2}} \int_{\Reals^d} |x|^{2s} e^{-x^2}dx ,\\
(C'_{d,s})^2=&\frac{1}{\pi^\frac{d}{2}} \int_{|x|<1} |x|^{2s} e^{-x^2}dx ,\\
r(a,t)=& \left( \frac{1+4a^2t^2}{a}\right)^\frac12 ,
\end{align*}
and let $\chi_{|x|<r(a,t)}$ be the characteristic function with support on $|x|<r(a,t)$. 

Given $A>0$ and $a>0$, let $\psi[A,a]$ be solutions to the \ls with initial data
\begin{align*}
u_0[A,a]=& A\left( \frac{a}{\pi}\right)^\frac{d}{4} e^{\frac{-ax^2}{2}} .
\end{align*}
These are given by 
\begin{align*}
\psi[A,a]= e^{it\Delta}u_0[A,a]=& A \left( \frac{a}{\pi}\right)^\frac{d}{4} \frac{1}{(1+2ait)^\frac{d}{2}} e^\frac{-ax^2}{2(1+2ait)} \\
\| e^{it\Delta}u_0[A,a] \|_{L^2_x} =& A\\
\| |x|^s e^{it\Delta}u_0[A,a] \|_{L^2_x} =& C_{d,s} A r(a,t)^s 
\sim A a^\frac{s}{2} t^s &\text{   for large $t$}
\end{align*}

Let $\psi_k=e^{it\Delta}u_0[A_k,a_k]$ with $A_j$ chosen, in terms of $A_k$ for $k<j$, to be sufficiently small so that 
\begin{align}
\label{eM4.1}
\sum_{i=k+1}^j A_i \leq \frac18 C'_{d,s} A_k
\end{align}
and with $a_j$ chosen sufficiently large so that 
\begin{align}
\label{eM4.2}
\sum_{k<j} A_k (1+C_{d,s}r(a_k,t)^s) \leq& \frac14 C'_{d,s} A_j a_j^\frac{s}{2}t^s ,
\end{align}
and let $\Psi= \sum \psi_k$. 

Given a fixed $t$, let $r_k= r(a_k,t)$ and $\chi_k=\chi_{|x|<r_k}$. The first condition, \eqref{eM4.1}, ensures that, for sufficiently large $k$, on a length scale of $|x|<r_k$, the function $\psi_k$ dominates all the later $\psi_j$ with $j>k$:
\begin{align*}
\| \chi_k |x|^s \psi_k\|_{L^2}
=C'_{d,s} A_k r_k^s > 8 \sum_{j>k} A_j r_k^s> 8r_k^s \sum_{j>k}\|\psi_j\|_{L^2}
>4\sum_{j>k} \| \chi_k |x|^s \psi_j\|_{L^2}
\end{align*}
The second condition, \eqref{eM4.2}, ensures that $A_j a_j^\frac{s}{2}$ grows at least exponentially. It also ensures that, for $a_j^\frac{1}{2}>t^{-2}$, in $\Xsigma$,$ \psi_j$ dominates all of the previous $\psi_k$ with $k<j$:
\begin{align*}
\| \chi_j |x|^2\psi_j\|_{L^2} 
=C'_{d,s}A_j r_j^s
>C'_{d,s}A_j a_j^\frac{s}{2}t^s >4\sum_{k<j}A_k (1+C_{d,s} r_k^s)
> 4\sum_{k<j} \|\psi_k\|_{\Xsigma}
\end{align*}

Since $\Psi$ is the sum of the $\psi_k$, and since, at a given time $t$, for sufficiently large $j$, on a length scale of $r_j$, $\psi_j$ dominates all the other $\psi_k$, the $\Xsigma$ norm of $\Psi(t)$ is bounded below by arbitrarily large numbers and must diverge:
\begin{align*}
\| \Psi(t)\|_{\Xsigma}
\geq \| \chi_j |x|^s \Psi\|_{L^2} \geq \| \chi_j |x|^s \psi_j\|_{L^2} - \sum_{k\not=j} \| \chi_j |x|^s \psi_k \|_{L^2} \geq \frac12 \| \chi_j |x|^s \psi_j\|_{L^2}
\geq C'_{d,s} A_j a_j^\frac{s}{2} t^s \rightarrow \infty .
\end{align*}
\end{proof}

A similar sequence of Gaussian initial data shows that the \nls equation is also ill-posed in $\Xsigma$. 

\begin{proposition}
\label{PM5.1}
\label{PXsigmaIllPosedness}
The \nls equation \eqref{NLS} is not well-posed in $\Xsigma$ for any $s>0$: the \nls evolution from time $0$ to time $t>0$, $S_{NLS}(t,0):u_0\mapsto u(t)$, is not continuous from $\Xsigma$ to $\Xsigma$. 
\end{proposition}
\begin{proof}It is sufficient to construct a sequence of $L^2$ solutions, $u^{[k]}$, for which $\|u^{[k]}_0\|_{\Xsigma}\rightarrow 0$, but for which, at any $t>0$, $\|u^{[k]}(t)\|_{\Xsigma}\rightarrow\infty$. This construction follows from the closeness of \nls and \ls evolutions for small initial data and from the existence of linear solutions with arbitrarily fast $\Xsigma$ growth. 

When the $L^\SBalancedExponent_{tx}$ norm is sufficiently small, the difference between the \ls and \nls evolutions is small. From Duhamel's principle and an extension of the local well-posedness theory, it is known that there is a $\delta'$ such that, if $\|u_0\|_{L^2}\leq\delta'$, then $u$ is defined for all time, and 
\begin{align*}
\| u-e^{it\Delta}u_0\|_{S^0}\leq& \| u\|_{L^\SBalancedExponent_{tx}}^{\frac4d+1} , \\
\|u\|_{L^\SBalancedExponent_{tx}}\leq& 2\|u_0\|_{L^2} .
\end{align*}

For the \ls solutions, the notation from the previous lemma will be used. In addition, $u^{[k]}$ will denote the \nls evolution of $u_0[A_k,a_k]$ with $A_k$ decreasing to zero and $a_k$ increasing to infinite, but with rates to be chosen. The index $k$ will be chosen sufficiently large so that, for $i>k$, 
\begin{align*}
A_i <& \delta'\\
A_i^\frac4d <& \frac{C'_{d,s}}{ 2^{\frac4d+2}} 
\end{align*}

As in Lemma \ref{LM3.1}, the $\Xsigma$ norm can be estimated by localizing on a length scale of $r_k$. 
\begin{align*}
\| u^{[k]}(t) \|_{\Xsigma}
\geq& \| \chi_k |x|^s u^{[k]} \|_{L^2}
\geq \| \chi_k |x|^s e^{it\Delta} u_0^{[k]} \|_{L^2} - \| \chi_k |x|^s (u^{[k]}-e^{-it\Delta} u_0^{[k]}) \|_{L^2}\\
\geq&  C'_{d,s} A_k r_k^s - r_k^s 2^\frac4d \| u_0^{[k]}\|_{L^2}^\frac4d 
\geq C'_{d,s} A_k r_k^s (1-\frac12)
\gtrsim A_k a_k^\frac{s}{2} t^s .
\end{align*}
If $A_k=a_k^\frac{-s}{4}$ and $a_k\rightarrow\infty$, then, for each $t>0$, 
\begin{align*}
\| u^{[k]}_0 \|_{\Xsigma}\sim& A_k\rightarrow 0 , \\
\| u^{[k]}(t) \|_{\Xsigma}\sim& A_ka_k^\frac{s}{2} t^s \rightarrow \infty .
\end{align*}
\end{proof} 

\subsection{The pseudoconformal transform and Strichartz norms}
The pseudoconformal transform is a symmetry of both the \ls equation and the pseudoconformal \nls equation \eqref{NLS} and is also an isometry on $L^2_x$ and the Strichartz admissible $L^qL^r$ spaces. 

\begin{lemma}
\label{T10.2}
\begin{enumerate}
\item If $0\leq t_1\leq t_2$, and $u:[t_1,t_2]\times\Reals^d\rightarrow\Complex$ satisfies the \nls equation \eqref{NLS}, then so does $\Conf[u]=v:[-t_1^{-1},-t_2^{-1}]\times\Reals^d\rightarrow\Complex$. 
\item If $u:\{t_1\}\times\Reals^d\rightarrow\Complex$, then $\|u(t_1)\|_{L^2}=\|v(-t_1^{-1} )\|_{L^2}$. 
\item If $u:[t_1,t_2]\times\Reals^d\rightarrow\Complex$ and $(q,r)$ is admissible, then $\|u\|_{L^q_{t\in[t_1,t_2]}L^r_x}=\|v\|_{L^q_{t\in[-t_1^{-1},-t_2^{-1}]}L^r_x}$.
\item Up to a reflection, the pseudoconformal transform is its own inverse: $\Conf[\Conf[u]](t,x)=u(t,-x)$
\end{enumerate}
\end{lemma}

These facts may be validated through explicit calculations.

\section{Global existence for initial data in $\Xsigma$}
\label{s2}

The goal is to prove global existence for initial data in $\Xsigma$ from the assumption that there is global well-posedness in $H^\RegExponent$. 

\JimsMacDoesntLikeFig{
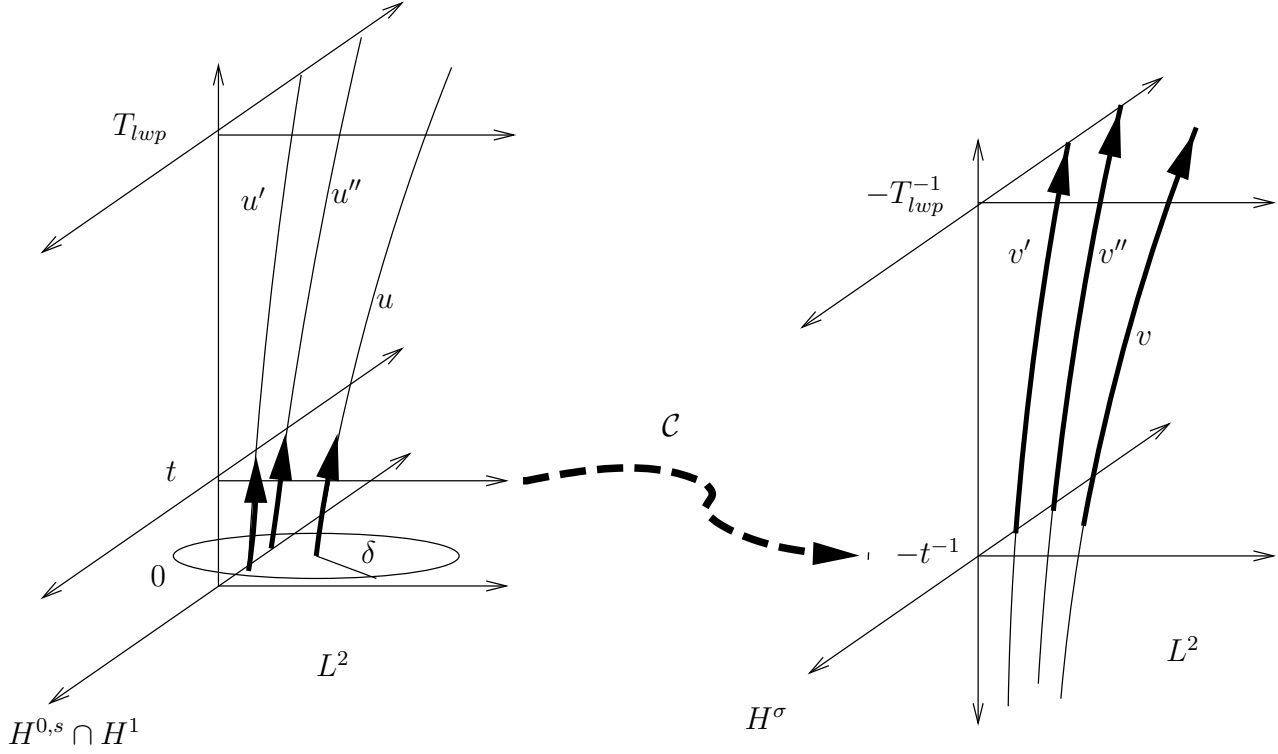
\begin{figure}
\input{NLSConformalPicture1.pstex_t}
\caption{The map $F_t:u_0\mapsto v(-T_{lwp}^{-1})$. 
$u_0$ is given in $\Xsigma\subset L^2$, and $u'_0$ and $u''_0$ in $\Xsigma\cap H^1$ are then chosen in a $\Xsigma$ neighborhood of $u$. }
\label{Fig1}
\end{figure}
}

Heuristically, initial data $u_0\in\Xsigma$ at $t_0=0$ can be transformed to initial data $v_0\in H^\RegExponent$ at $\transT_0=-\infty$. Under the $H^\RegExponent$ global well-posedness hypothesis, $v$ can then be defined for all time, and $u$ can be defined for all time by the inverse pseudoconformal transform. To make this heuristic rigorous, $u_0$ can be evolved to $u(T_{lwp})$ and then pseudoconformally transformed to $v$. Following this, it is sufficient to show that $v(-T_{lwp}^{-1})$ is in $H^\RegExponent$ to apply the $H^\RegExponent$ global well-posedness hypothesis. 

In terms of the \nls evolution map, $S_{NLS}(t_2,t_1)$, which was introduced earlier, the map $u_0\mapsto v(-T_{lwp}^{-1})$ is
\begin{equation*}
F=\Conf\circ S_{NLS}(T_{lwp},0): L^2\rightarrow L^2 .
\end{equation*}

Since the pseudoconformal transform commutes with the \nls evolution, this map can also be constructed in a different way, which is illustrated in Figure \ref{Fig1}. The solutions are first allowed to evolve under $S_{NLS}(t,0)$ to time $t$ (dark in left diagram), then pseudoconformally transformed from data at time $t$ to data at transformed time $-t^{-1}$ (dashed arrow from left diagram to right diagram), and finally allowed to evolve under $S_{NLS}(-T_{lwp}^{-1},-t^{-1})$ to transformed time $-T_{lwp}^{-1}$ (dark in right diagram). This construction is represented by 
\begin{equation*}
F=F_t=S_{NLS}(-T_{lwp}^{-1},-t^{-1})\circ\Conf\circ S_{NLS}(t,0):L^2\rightarrow L^2 .
\end{equation*}
By the $L^2$ local well-posedness Theorem \ref{T6.1} and the properties of the pseudoconformal transform, in a $L^2$ neighborhood of $u_0$, $F$ is continuous with respect to the $L^2$ norm. 

To prove Proposition \ref{P3.1}, it is sufficient to show that $F$ can be restricted to $F:\Xsigma\rightarrow H^\RegExponent$. This is done by initially restricting to regularized data in $\Xsigma\cap H^1$, showing that each of the three steps of $F_t$ is continuous with respect to the regularized data, and then removing the regularization. $\Xsigma\cap H^1$ is a useful auxiliary space because it is preserved by the \nls evolution and pseudoconformally transforms to $H^\RegExponent$. 

Because the $H^\RegExponent$ local well-posedness Theorem \ref{T7.1} uses $L^2$ Strichartz norms to control the divergence of nearby solutions, if $u'$ and $u''$ start near $u$ in $L^2$, their separation in $\Xsigma\cap H^1$ can not grow by more than a constant factor. Similarly, if $v'$ and $v''$ start near $v$ in $L^2$, then their $H^\RegExponent$ separation can not increase by more than a constant factor. Thus the divergence of the approximators, from each other, is controlled if they start in a sufficiently small $L^2$ neighborhood of $u_0$. 

This $L^2$ neighborhood of $u_0$ is illustrated by the oval in the left diagram in Figure \ref{Fig1}. It is taken to be a ball of radius $\delta$, and this is the $\delta$ which appears in the following subsections. The value of $\delta$ is dictated by the $H^\RegExponent$ local well-posedness Theorem \ref{T7.1}. 

Since $F=F_t$ is independent of $t$, it is possible to take the infimum in $t$ of the $H^\RegExponent$ norm estimates for $v'(-T_{lwp}^{-1})$ and $v''(-T_{lwp}^{-1})$. This eliminates the dependence on the $H^1$ regularization and corresponds to the original heuristic of transforming $\Xsigma$ data at time $t_0=0$ to $H^\RegExponent$ data at $\transT_0=-\infty$. As a result $F:\Xsigma\cap H^1\rightarrow H^\RegExponent$ is continuous with respect to the $\Xsigma$ norm alone and extends uniquely to $F:\Xsigma\rightarrow H^\RegExponent$ in a neighborhood of $u_0$. As a result, and contrary to the image in Figure \ref{Fig1}, $v(-T_{lwp}^{-1})$ is in $H^\RegExponent$ not merely $L^2$. 

In Subsection \ref{ss2.2}, it will be shown that the \nls evolution is continuous in $\Xsigma\cap H^1$. In Subsection \ref{ss2.3}, real interpolation will be used to show that the pseudoconformal transform takes $\Xsigma\cap H^1$ to $H^\RegExponent$ with a $t$ dependent coefficient on the $H^1$ part of the norm. In Subsection \ref{ss2.4}, the $H^\RegExponent$ local well-posedness theorem will be used to show that $H^\RegExponent$ data evolves continuously in $H^\RegExponent$ from transformed time $-t^{-1}$ to transformed time $-T_{lwp}^{-1}$. In Subsection \ref{ss2.5}, the infimum in $t$ will be taken to eliminate the $H^1$ dependence. 

\subsection{$S_{NLS}(t,0):\Xsigma\cap H^1\rightarrow \Xsigma\cap H^1$}
\label{ss2.2}

By the local well-posedness theory, the \nls evolution takes $H^1\rightarrow H^1$, at least up to the $L^2$ local well-posedness time, $T_{lwp}$. For the linear evolution, a simple commutator calculation shows that the $H^1$ norm controls the growth of the $\Xsigma$ norm. The extra terms arising in the nonlinear evolution cancel, so that the same result holds. 

\begin{lemma}
\label{L13.1}
For $u$ a solution to the \nls equation \eqref{NLS} with initial
data $u_0\in \Xsigma$ with $0 \leq s \leq 1$ and local existence time $T_{lwp}$, there is a $\delta$ \deltacomesfromphrase such that
\begin{enumerate}
\item if $u'$ solves the \nls equation with initial data $u'_0\in\Xsigma\cap H^1$ and with $\|u_0-u'_0\|_{L^2}<\delta$, then $\forall~ t\in[0,T_{lwp}]$
\begin{align*}
\|u'(t)\|_{H^1}\leq& 2\|u'_0\|_{H^1} , \\
\|u'(t)\|_{\Xsigma}\leq& \|u'_0\|_{\Xsigma} + 2t\|u'_0\|_{H^1}.
\end{align*}
\item if $u'$ and $u''$ solve the \nls equation with
initial data $u'_0\in\Xsigma\cap H^1$ and $u''_0 \in\Xsigma\cap H^1$
respectively and with $\|u_0-u'_0\|_{L^2}+ \|u_0-u''_0\|_{L^2} + \|u'_0-u''_0\|_{L^2}<\delta$, then $\forall~ t\in[0,T_{lwp}]$, 
\begin{align*}
\|u'(t)-u''(t)\|_{H^1}\leq& 2\|u'_0-u''_0\|_{H^1}\\
\|u'(t)-u''(t)\|_{\Xsigma}\leq& \|u'_0-u''_0\|_{\Xsigma} + 2t\|u'_0\|_{H^1}+ 2t\|u''_0\|_{H^1} .
\end{align*}
\end{enumerate}
\end{lemma}
\begin{remark}
The solutions $u'$ and $u''$ are thought of as perturbations from $u$. The original solution $u$ appears in the statement of this theorem because it provides both the local well-posedness interval, $[0,T_{lwp}]$, on which we wish to control $u'$ and $u''$, and the $L^\SBalancedExponent_{tx}$ estimates with which we can achieve this control. 
\end{remark}
\begin{proof}
From the local well-posedness Theorem \ref{T6.1}, $T_{lwp}$ can be chosen small enough so that $\|u\|_{L^\SBalancedExponent_{tx}([0,T_{lwp}])}$ is less than half of $\delta_3$ from the $H^s$ local well-posedness Theorem \ref{T7.1}. In this case, by the $L^2$ local well-posedness Theorem \ref{T6.1}, $\|u-u'\|_{L^\SBalancedExponent_{tx}([0,T_{lwp}])}\leq 2\|u_0-u'_0\|_{L^2}<2\delta$ and $\|u'\|_{L^\SBalancedExponent_{tx}([0,T_{lwp}])}<2\delta+\frac12\delta_3$. 

The function $u'$ can now be taken as the solution to estimate and from which $u''$ is a perturbation. If $2\delta+\frac12\delta_3\leq\delta_3$, then the $H^\RegExponent$ local well-posedness Theorem \ref{T7.1} provides the estimates on the growth and separation in $H^1$. 

Differentiation in time and the Cauchy-Schwartz estimate gives the growth of the weighted norms. 
\begin{align*}
\Dt \|u'(t)\|_{\Xsigma}^2
=& \Dt \langle u', \Jax^{2\RegExponent} u'\rangle . \\
2\|u'(t)\|_{\Xsigma}\Dt \|u'(t)\|_{\Xsigma}
=& \langle i\Delta u' - i\lambda |u'|^{\frac4d}u',\Jax^{2\RegExponent}u' \rangle +\langle u',\Jax^{2\RegExponent}(i\Delta u' -i\lambda |u'|^{\frac4d}u') \rangle\\
=& -i\langle \Delta u',\Jax^{2\RegExponent}u'\rangle +i\langle u',\Jax^{2\RegExponent}\Delta u'\rangle\\
& + i\lambda\langle |u'|^{\frac4d}u',\Jax^{2\RegExponent} u'\rangle - i\lambda\langle u',\Jax^{2\RegExponent}|u'|^{\frac4d}u'\rangle \\
=& \langle u', i[-\Delta,\Jax^{2\RegExponent}]u'\rangle . \\
|2\|u'(t)\|_{\Xsigma}\Dt \|u'(t)\|_{\Xsigma}|
\leq& | \langle u', -4i\RegExponent\Jax^{2\RegExponent-2 }x D u'  -2i\RegExponent\Jax^{2\RegExponent-2}((2\RegExponent-2)\frac{x^2}{\Jax^2} + 1)u'\rangle | \\
\leq& C\|\Jax^{2\RegExponent-1}u'\|_{L^2} \|D_x u'\|_{L^2} + C\|\Jax^{\RegExponent-1} u'\|_{L^2}^2\\
\leq& C\|u'\|_{\Xsigma} \|D_x u'\|_{L^2} + C\|u'\|_{\Xsigma}\|u'\|_{L^2} . \\
\Dt\|\Jax u'\|_{\Xsigma}
\leq& C \|D_x u'\|_{L^2} +C \|u'\|_{L^2}\\
\leq& C \|u'\|_{H^1} .
\end{align*}
This proves that 
\begin{align*}
\|u'(t)\|_{\Xsigma} \leq \|u'_0\|_{\Xsigma} + Ct\|u'_0\|_{H^1} .
\end{align*}
A similar calculation shows 
\begin{align*}
\|u'(t)-u''(t)\|_{\Xsigma}
\leq& \|u'_0-u''_0\|_{\Xsigma} + \|u'(t)-u'_0\|_{\Xsigma} + \|u''(t)-u''_0\|_{\Xsigma}\\
\leq& \|u'_0-u''_0\|_{\Xsigma} + Ct(\|u'_0\|_{H^1}+\|u''_0\|_{H^1}) .
\end{align*}
\end{proof}

\subsection{$\Conf:\Xsigma\cap H^1\rightarrow H^\RegExponent$}
\label{ss2.3}

In this section, it is shown that the pseudoconformal transform takes a function $u(t)\in\Xsigma\cap H^1$ to $v(-t^{-1})\in H^\RegExponent$. This is done by interpolation between $L^2$ and $\Xone\cap H^1$ using the $K$ method of real interpolation. 

To begin, the arguments for $L^2$ and $\Xone\cap H^1$ are presented. The $L^2$ result is part of Theorem \ref{T10.2}. The $\Xone\cap H^1$ result leads to equation \eqref{3.1}, which was stated in the introduction. 

\begin{lemma}
\label{L14.1}
If $u:\{t\}\times\Reals^d\rightarrow\Complex$ and $v=\Conf[u]$, then
\begin{align}
\label{14.1}
\|v(-t^{-1})\|_{L^2}=&\|u(t)\|_{L^2} , \\
\label{14.2}
\|v(-t^{-1})\|_{H^1}\leq&\|u(t)\|_{\Xone}+t\|u(t)\|_{H^1}
\end{align}
where it is understood that the norm on the right is infinite (and the inequality trivial) if $u(t)$ does not belong to the appropriate space. 
\end{lemma}
\begin{proof}
Equation \eqref{14.1} is statement 2 of Lemma \ref{T10.2}.

Inequality \eqref{14.2} follows by direct computation with a change of variables. This computation is simplified by recalling the notation $\transT=-t^{-1}$ and $ty=x$. 
\begin{align*}
\|v(-t^{-1})\|_{H^1}^2
=&\int |\nabla_y v(-t^{-1},y)|^2 dy\\
=&\int |\nabla_y(t^\frac{d}{2} e^\frac{it y^2}{4} u(t,ty))|^2 dy\\
=&\int t^d |\frac12 ity e^\frac{ity^2}{4}u(t,ty)+e^\frac{ity^2}{4}\nabla_yu(t,ty)|^2 dy\\
\leq&\int |ty u(t,ty)|^2 t^d dy+\int 2|\nabla_y u(t,ty)|^2 t^d dy\\
\leq&\int |xu(t,x)|^2dx + 2\int t^2|\nabla_xu(t,x)|^2 dx.
\end{align*}
\end{proof}

The $K$ method of real interpolation is now summarized from \cite{BerghLofstrom}. The $\RegExponent$ interpolation norm of $a\in A_0+A_1$ is defined by the following, if this norm is finite, 
\begin{align}
\label{15.1}
K(\lambda,a;A_0,A_1)=&\inf_{a=a_0+a_1}(\|a_0\|_{A_0}^2+\lambda^2\|a_1\|_{A_1})^2\\
\label{15.2}
\|a\|_{\RegExponent, (A_0,A_1); K}^2=&\Phi_{\RegExponent,2}(K(\lambda,a;A_0,A_1))\\
=&\int \lambda^{-2\RegExponent-1}K(\lambda,a;A_0,A_1)^{2} d\lambda. \nonumber
\end{align}
Since only the $K$ method of interpolation will be introduced, the $K$ index in the norm will be omitted $\|a\|_{\RegExponent,(A_0,A_1)}=\|a\|_{\RegExponent, (A_0,A_1); K}$. If $a\in A_0\cap A_1$, then 
\begin{align*}
\|a\|_{\RegExponent,(A_0,A_1)}\leq \|a\|_{A_0}^{1-\RegExponent} \|a\|_{A_1}^\RegExponent .
\end{align*}
The interpolation space $(A_0,A_1)_\RegExponent$ is defined as the set of $a\in A_0+A_1$ for which $\|a\|_{\RegExponent,(A_0,A_1)}$ is finite. There are some technical issues, but since only spaces $A_0$ and $A_1$ which are subsets of $L^2$ will be considered, $(A_0,A_1)_\RegExponent$ will be well-defined, a Banach space, and the closure of $A_0\cap A_1$. 

The $K$ method of real interpolation is an exact interpolation method of exponent $\RegExponent$ \cite{BerghLofstrom}: if
\begin{align*}
T: A_0&\rightarrow B_0 & \|T\|_{A_0\rightarrow B_0}\leq M_0\\
T: A_1&\rightarrow B_1 & \|T\|_{A_1\rightarrow B_1}\leq M_1 ,
\end{align*}
then
\begin{align}
T: (A_0,A_1)_\RegExponent\rightarrow (B_0,B_1)_\RegExponent \nonumber \\
\label{15.3}
\|T\|_{(A_0,A_1)_\RegExponent\rightarrow (B_0,B_1)_\RegExponent} \leq
M_0^{1-\RegExponent} M_1^\RegExponent .
\end{align}

It is known that 
\begin{align*}
(L^2,H^1)_\RegExponent=&H^\RegExponent\\
(L^2,\Xone)_\RegExponent=&\Xsigma .
\end{align*}

By interpolating the results of Lemma \ref{L14.1}, it follows that
\begin{align*}
\Conf:(L^2,\Xone\cap H^1)_\RegExponent\rightarrow H^\RegExponent .
\end{align*}
Unfortunately, because of the $\inf$ in \eqref{15.1}, it is not clear that $\|a\|_{\RegExponent,(L^2,\Xone\cap H^1)} = \|a\|_{\RegExponent,(L^2,\Xone)} + \|a\|_{\RegExponent,(L^2, H^1)} = \|a\|_{\Xsigma} + \|a\|_{H^\RegExponent}$; although, we expect this is true. We will instead prove the simpler result that 
\begin{equation*}
\|a\|_{\RegExponent,(L^2,\Xone\cap H^1)} \lesssim \|a\|_{\Xsigma} + \|a\|_{H^1} .
\end{equation*}
 To ensure that the $t$ dependent coefficients only appear on the $H^1$ norm, the $t$ dependence is kept in the interpolation calculations rather than being estimated by \eqref{15.3}. 

\begin{lemma}
\label{L16.1}
If $u:\{t\}\times\Reals^d\rightarrow\Complex$, $u(t)\in\Xsigma\cap H^1$, and $v=\Conf[u]$, then
\begin{equation*}
\|v(-t^{-1})\|_{H^\RegExponent}\leq C_{1,\RegExponent} \|u\|_{\Xsigma} + C_{2,\RegExponent}t \|u\|_{H^1}
\end{equation*}
where the constants $C_{1,\RegExponent}$ and $C_{2,\RegExponent}$ depend only on
$\RegExponent$ and $d$. 
\end{lemma}
\begin{proof}
Let $u:\{t\}\times\Reals^d\rightarrow\Complex$ with $u(t)$ in Schwartz class. Using the $K$-method, it will be shown that the $H^\RegExponent$ norm is dominated by the $\Xsigma$ and $H^1$ norms. The $K$-method of interpolation involves taking an infimum over all possible decompositions of $u$. This infimum is dominated by any particular choice of decomposition. The decomposition which is optimal for balancing $L^2$ with $\Xone$ will be used. This will give the $\Xsigma$ part of the estimate. There is no reduction in the regularity required for the estimate, since this decomposition ignores the $H^1$ term. 

\begin{align*}
\|v(-t^{-1})\|_{H^\RegExponent}
=& \|\Conf[u](-t^{-1})\|_{H^\RegExponent}\\
=&\int \lambda^{-2\RegExponent-1} K(\lambda,\Conf[u];L^2,H^1)^2 d\lambda
\end{align*}

Using Lemma \ref{L14.1}, $K(\lambda,\Conf[u])$ can be estimated in terms of the $L^2$, $\Xone$, and $H^1$ norms of $u$. (Note that in this proof, $u_0$ refers to part of the interpolation decomposition in \eqref{15.1}, not the initial data.) 
\begin{align*}
K(\lambda,\Conf[u])
=&\inf_{u=u_0+u_1}(\|\Conf[u_0]\|_{L^2}^2 +\lambda^2 \|\Conf[u_1]\|_{H^1}^2)^\frac12\\
\leq& \inf_{u=u_0+u_1}(\|u_0\|_{L^2} +
\lambda^2\|u_1\|_{\Xone}^2+\lambda^2\|u_1\|_{H^1}^2)^\frac12 .
\end{align*}
In the proof that $(L^2,\Xone)_{\RegExponent}=\Xsigma$ \cite{BerghLofstrom}, it is shown that the optimal decomposition for $\|u_0\|_{L^2} + \lambda^2\|u_1\|_{\Xone}^2$ is
\begin{align*}
u_0=& \frac{\lambda^2(1+x^2)}{1+\lambda^2(1+x^2)} u\\
u_1=&\frac{1}{1+\lambda^2(1+x^2)} u .
\end{align*}
This decomposition will be used for $\lambda<1$ to bound $K(\lambda,\Conf[u])$ from above. 
\begin{align*}
K(\lambda,\Conf[u])^2
\leq& \|u_0\|_{L^2}^2 +\lambda^2\|u_1\|_{\Xone}^2 +\lambda^2\|u_1\|_{H^1}^2\\
\leq& \int ( \frac{\lambda^4(1+x^2)^2}{(1+\lambda^2(1+x^2))^2} + \frac{1}{(1+\lambda^2(1+x^2))^2}\lambda^2(1+x^2) )|u|^2 d^2x\\
&+\lambda^2 t^2 \int |\nabla (\frac{1}{1+\lambda^2(1+x^2)} u)|^2 dx\\
\leq& \int \frac{\lambda^2(1+x^2)}{1+\lambda^2(1+x^2)} |u|^2 dx \\
&+\lambda^2 t^2 \int |\nabla (\frac{1}{1+\lambda^2(1+x^2)} u)|^2 dx .
\end{align*}
This decomposition can be used to bound the $H^\RegExponent$ norm of $v=\Conf[u]$ for $\lambda<1$. The decomposition $u_0=u$ and $u_1=0$ will be used for $\lambda\geq1$. 
\begin{align*}
\| v(-t^{-1})\|_{H^\RegExponent}^2
\leq& \int_0^1 \lambda^{-2\RegExponent-1}K(\lambda,\Conf[u]) d\lambda + \int_1^\infty \lambda^{-2\RegExponent-1}K(\lambda,\Conf[u]) d\lambda\\
\leq& \int_0^1 \lambda^{-2\RegExponent-1} \int_{\Reals^2} \frac{\lambda^2(1+x^2)}{1+\lambda^2(1+x^2)} |u|^2 dx d\lambda \\
&+\int_0^1 \lambda^{-2\RegExponent-1} \int \lambda^2 t^2|\nabla_x (\frac{1}{1+\lambda^2(1+x^2)} u)|^2 dx d\lambda\\
&+\int_1^\infty \lambda^{-2\RegExponent-1} \|u_0\|_{L^2}^2 d\lambda\\
\leq& \int_0^1 \lambda^{-2\RegExponent-1} \int_{\Reals^2} \frac{\lambda^2(1+x^2)}{1+\lambda^2(1+x^2)} |u|^2 dx d\lambda \\
&+ \int_0^1 \lambda^{-2\RegExponent+1} t^2 \int \lambda |u|^2+|\nabla_x u|^2 dx d\lambda\\
&+ C_\RegExponent \|u\|_{L^2}^2 .
\end{align*}
At this stage, the first term is evaluated by the substitution $\lambda'^2=(1+x^2)\lambda^2$ and Fubini's theorem. The other two pieces are estimated by direct integration and estimated using the assumption $t<1$. 
\begin{align*}
\| v(-t^{-1})\|_{H^\RegExponent}^2
\leq& \int_{\Reals^2} (1+x^2)^\RegExponent |u|^2 \int_0^1 \lambda'^{-2\RegExponent-1}\frac{\lambda'^2}{1+\lambda'^2} d\lambda' dx\\
&+ Ct^2 \|u\|_{H^1}^2 + C'(1+t^2)\|u\|_{L^2}^2\\
\leq& C_{1,\RegExponent} \|u\|_{\Xsigma}^2 + C_{2,\RegExponent}t^2 \|u\|_{H^1}^2 .
\end{align*}

\end{proof}

\subsection{$S_{NLS}(-T_{LWP}^{-1},-t^{-1}):H^\RegExponent\rightarrow H^\RegExponent$}
\label{ss2.4}
Since the pseudoconformal transform preserves Strichartz admissible norms, $\|v\|_{L^\SBalancedExponent_{tx}([-t^{-1},-T_{lwp}^{-1}])}$ is controlled by $\|u\|_{L^\SBalancedExponent_{tx}([0,T_{lwp}])}$ and can be taken to be small. Since the growth of the $H^\RegExponent$ norm under the \nls evolution is controlled by $L^\SBalancedExponent_{tx}$, $\|v(-T_{lwp}^{-1})\|_{H^\RegExponent}$ can be controlled by $\|v(-t^{-1})\|_{H^\RegExponent}$ for any $t\in[0,T_{lwp}]$. The divergence of nearby solutions can be similarly controlled. From this, if $u'_0$ and $u''_0$ are $\Xsigma\cap H^1$ approximations of $u_0\in\Xsigma$, then the $H^\RegExponent$ distance between $v'$ and $v''$ at transformed time $-T_{lwp}^{-1}$ is controlled by the distance between them at earlier transformed times. 

\begin{lemma}
\label{L18.1}
For $u$ a solution to the \nls equation \eqref{NLS} with initial data $u_0\in L^2$, there is a local well-posedness time $T_{lwp}$ and a $\delta$ \deltacomesfromphrase such that:
\begin{enumerate}
\item if $u'$ is a solution to the \nls equation on $[0,T_{lwp}]$ with
  initial data $u'_0 \in \Xsigma \cap H^1 $,
  with $\|u_0-u'_0\|_{L^2}<\delta$ and with a $t\in[0,T_{lwp}]$ such that $v'(-t^{-1})\in H^\RegExponent$, then 
\begin{align*}
\|v'(-T_{lwp}^{-1})\|_{H^\RegExponent} <2\|v'(-t^{-1})\|_{H^\RegExponent} ;
\end{align*}
\item if $u'$ and $u''$ are solutions to the \nls equation on $[0,T_{lwp}]$ with initial data $u'_0\in L^2$ and $u''_0\in L^2$ respectively, with $\|u_0-u'_0\|_{L^2}<\delta$ and $\|u_0-u''_0\|_{L^2}<\delta$, and with a $t\in[0,T_{lwp}]$ such that $v'(-t^{-1})\in H^\RegExponent$ and $v''(-t^{-1})\in H^\RegExponent$, then
\begin{align*}
\|v'(-T_{lwp}^{-1})-v''(-T_{lwp}^{-1})\|_{H^\RegExponent} <2\|v'(-t^{-1}
)-v''(-t^{-1} )\|_{H^\RegExponent} .
\end{align*}
\end{enumerate}
\end{lemma}
\begin{remark}
As in Lemma \ref{L13.1}, $u$ appears in this theorem because we want to control the size of $L^\SBalancedExponent_{tx}$ norms for the approximators, $u'$ and $u''$. 
\end{remark}
\begin{proof}
From the $L^2$ local well-posedness Theorem \ref{T6.1}, $T_{lwp}$ can be chosen small enough so that $\|u\|_{L^\SBalancedExponent_{tx}([0,T_{lwp}])}$ is less than half of $\delta_3$ from the $H^\RegExponent$ local well-posedness Theorem \ref{T7.1}. In this case, by the $L^2$ local well-posedness Theorem \ref{T6.1}, $\|u-u'\|_{L^\SBalancedExponent_{tx}([0,T_{lwp}])}\leq 2\|u_0-u'_0\|_{L^2}<2\delta$ and $\|u'\|_{L^\SBalancedExponent_{tx}([0,T_{lwp}])}<2\delta+\frac12\delta_3$. 

The function $u'$ or $v'$ can now be taken as the solution to estimate. Since the pseudoconformal transformation preserves the Strichartz admissible norms, $\|v'\|_{L^\SBalancedExponent_{tx}([-t^{-1},T_{lwp}^{-1}])} \leq \|v'\|_{L^\SBalancedExponent_{tx}((-\infty,-T_{lwp}^{-1}])} =\|u'\|_{L^\SBalancedExponent_{tx}([0,T_{lwp}])}<2\delta+\frac12\delta_3$. Thus, if $2\delta+\frac12\delta_3\leq\delta_3$, then, by the $H^\RegExponent$ local well-posedness Theorem \ref{T7.1}, 
\begin{align*}
\|v'\|_{S^\RegExponent([-t^{-1},-T_{lwp}^{-1}])} <2\|v'(-t_{-1})\|_{H^\RegExponent}.
\end{align*}

If $\|u_0-u'_0\|_{L^2}\leq\delta$ and $\|u_0-u''_0\|_{L^2}\leq\delta$, then $\|u'_0-u''_0\|_{L^2}\leq2\delta$. Again, by the $H^\RegExponent$ local well-posedness Theorem \ref{T7.1}, if, in addition to the previous conditions, $2\delta\leq\delta_4$, then
\begin{align*}
\|v'(-T_{lwp}^{-1})-v''(-T_{lwp}^{-1})\|_{H^\RegExponent}\leq \|v'-v''\|_{S^\RegExponent([-t^{-1},-T_{lwp}^{-1}])}\leq 2\|v'(-t^{-1})-v''(-t^{-1})\|_{H^\RegExponent}.
\end{align*}
\end{proof}

\subsection{$F:\Xsigma\cap H^1\rightarrow H^\RegExponent$}
\label{ss2.5}

The results from Sections \ref{ss2.2}, \ref{ss2.3}, and \ref{ss2.4} are now combined to show that $F=F_t=S_{NLS}(-T_{lwp}^{-1},-t^{-1})\circ\Conf\circ S_{NLS}(t,0)$ takes $\Xsigma\cap H^1$ to $H^\RegExponent$. Furthermore, an explicit $t$ dependence on the $H^\RegExponent$ norm will be found in Lemma \ref{L20.1} and then removed in Proposition \ref{P22.1} to show that $F$ takes a $\Xsigma$ neighborhood of $u_0$ into $H^\RegExponent$. 

We will first show that for $u_0\in\Xsigma$, if we restrict attention to initial data which is both in a $\Xsigma$ neighborhood, $N$, of $u_0$ and in $H^1$, then $F$ maps this initial data in $N\cap H^1$ to $H^\RegExponent$. 

\begin{lemma}
\label{L20.1}
For $u$ a solution to the \nls equation \eqref{NLS} with initial data $u_0\in \Xsigma$, then there is a local well-posedness time $T_{lwp}$ and a $\delta$ \deltacomesfromphrase such that: 
\begin{enumerate} 
\item if $u'$ is a solution to the \nls equation with initial data
  $u'_0\in\Xsigma\cap H^1$ and $\|u_0-u'_0\|_{L^2}<\delta$, then for
  all $t \in (0, T_{lwp}]$ we have 
\begin{align}
\label{20.1a}
\|v'(-T_{lwp}^{-1})\|_{H^\RegExponent}\leq C_\RegExponent \|u'_0\|_{\Xsigma} +C_\RegExponent' t\|u'_0\|_{H^1} ;
\end{align}
\item if $u'$ and $u''$ are solutions to the \nls equation with
  initial data $u'_0\in\Xsigma\cap H^1$ and $u''_0\in\Xsigma\cap H^1$
  respectively and with $\|u_0-u'_0\|_{L^2}<\delta$ and
  $\|u_0-u''_0\|_{L^2}<\delta$, then for
  all $t \in (0, T_{lwp}]$ we have 
\begin{align}
\label{20.1b}
\|v'(-T_{lwp}^{-1})-v''(-T_{lwp}^{-1})\|_{H^\RegExponent}\leq C_\RegExponent \|u'_0-u''_0\|_{\Xsigma} +C_\RegExponent' t\|u'_0\|_{H^1}+C_\RegExponent' t\|u''_0\|_{H^1} .
\end{align}
\end{enumerate} 

In other words, there is an open set $N\in\Xsigma$ containing $u_0$ for which
\begin{align*}
F_t=S_{NLS}(-T_{lwp}^{-1},-t^{-1})\circ\Conf\circ S_{NLS}(t,0):N\cap H^1\rightarrow H^\RegExponent
\end{align*}
and $F_t$ is continuous with respect to the $\Xsigma\cap H^1$ topology. 
\end{lemma}
\begin{proof} 
Conditions on $\delta$ and $T_{lwp}$ will be found. To begin, assume $T_{lwp}<1$. 

Since $u_0\in\Xsigma$ and $u_0'$ and $u''_0$ are $\Xsigma\cap H^1$ approximations, by Lemma \ref{L13.1}, if $\delta$ is less than the $\delta$ in Lemma \ref{L13.1}, then, for $t\in[0,T_{lwp}]$, 
\begin{align*}
\|u'(t)\|_{H^1}\leq& 2\|u'_0\|_{H^1} ,  \\
\|u'(t)\|_{\Xsigma}\leq& \|u'_0\|_{\Xsigma} +2t \|u'_0\|_{H^1} , \\
\|u'(t)-u''(t)\|_{H^1}\leq& 2\|u'_0\|_{H^1}+2\|u''_0\|_{H^1}  ,\\
\|u'(t)-u''(t)\|_{\Xsigma}\leq& \|u'_0-u''_0\|_{\Xsigma} 
+2t \|u'_0\|_{H^1}+2t \|u''_0\|_{H^1} .
\end{align*}

By Lemma \ref{L16.1}, the linearity of the pseudoconformal transform
and the triangle inequality,
\begin{align*}
\|v'(-t^{-1})\|_{H^\RegExponent}
\leq& C\|u'_0\|_{\Xsigma} + tC'\|u'_0\|_{H^1}\\
\|v'(-t^{-1})-v''(-t^{-1})\|_{H^\RegExponent}
\leq& C\|u'_0-u''_0\|_{\Xsigma} + tC'\|u'_0\|_{H^1}+
tC'\|u''_0\|_{H^1} .
\end{align*}

Since $u_0$ is a solution with initial data in $\Xsigma$, $u'_0$ and $u''_0$ are $\Xsigma\cap H^1$ approximations, and $v'$ and $v''$ are in $H^\RegExponent$ at transformed time $-t^{-1}$, if $T_{lwp}$ and $\delta$ are less than the corresponding values in Lemma \ref{L18.1}, then 
\begin{align*}
\|v'(-T_{lwp}^{-1})\|_{H^\RegExponent} 
\leq& C\|u'_0\|_{\Xsigma} + tC'\|u'_0\|_{H^1} , \\
\|v'(-T_{lwp}^{-1})-v''(-T_{lwp}^{-1})\|_{H^\RegExponent}
\leq& C\|u'_0-u''_0\|_{\Xsigma} + tC'\|u'_0\|_{H^1}+ tC'\|u''_0\|_{H^1} .
\end{align*}

Since $F_t: u_0\mapsto v(-T_{lwp}^{-1})$ and the set $\|u_0-u'_0\|_{L^2}<\delta$ is open in $\Xsigma$, this set is the $N$ given in the statement of the theorem. By \eqref{20.1b}, $F_t$ is continuous from $\Xsigma\cap H^1$ to $H^\RegExponent$. 
\end{proof}

The infimum in $t$ can be taken when estimating the $H^\RegExponent$
norm of $F_t(u_0)=v'(-T_{lwp}^{-1})$. Since $F=F_t: u_0\mapsto
v(-T_{lwp}^{-1})$ is independent of $t$, this eliminates the $H^1$
dependence. Eliminating the $H^1$ dependence shows that $F$ is
continuous from $\Xsigma$ to $H^\RegExponent$. If $u_0$ is
approximated in $\Xsigma$ by a sequence of regularized initial data
$u^{[i]}_0$ in $\Xsigma \cap H^1$, then the corresponding $v^{[i]}(-T_{lwp}^{-1})$ must converge in $H^\RegExponent$ to $v(-T_{lwp}^{-1})$. This proves that $v(-T_{lwp}^{-1})$ is in $H^\RegExponent$. Under the assumption of $H^\RegExponent$ global existence, $v$ and $u$ can be extended globally. 

\begin{proposition}
\label{P22.1}
Assume that the \nls equation \eqref{NLS} is globally well-posed in $H^\RegExponent$ (with the additional hypothesis that the initial data has $L^2$ norm bounded by $\| Q \|_{L^2}$ in the focusing case). 

If $u_0\in\Xsigma$ (and $\|u_0\|_{L^2}<\| Q \|_{L^2}$ in the focusing case), then there  is a function $u:\Reals\times\Reals^d\rightarrow\Complex$ which solves the \nls equation \eqref{NLS} for all time. 

From the local well-posedness theory, this solution remains in $L^2$,
at each time, with constant norm and is the unique solution in the Strichartz space $S^0$. 
\end{proposition}
\begin{proof}
For $u_0\in\Xsigma$ at $t_0=0$, by Lemma \ref{L20.1}, there is a solution $u$ with local existence time $T_{lwp}$. If $u'_0$ and $u''_0$ are in $\Xsigma\cap H^1$ and sufficiently close to $u_0$ in $\Xsigma$, then, since the map $F:u_0\mapsto v(-T_{lwp}^{-1})$ is independent of $t$, it is possible to apply the infimum in $t$ to \eqref{20.1a} and \eqref{20.1b} and obtain
\begin{align*}
\|v'(-T_{lwp}^{-1})\|_{H^\RegExponent} \leq& C_{\RegExponent} \|u'_0\|_{\Xsigma}\\
\|v'(-T_{lwp}^{-1})-v''(-T_{lwp}^{-1})\|_{H^\RegExponent} \leq& C_{\RegExponent} \|u'_0-u''_0\|_{\Xsigma} .
\end{align*}
Since $\Xsigma\cap H^1$ is dense in $\Xsigma$, and $F$ is continuous with respect to the $\Xsigma$ norm, if a sequence $u_0^{[i]}\in\Xsigma\cap H^1$ is chosen to converge to $u_0$ in $\Xsigma$, the $v^{[i]}(-T_{lwp}^{-1})$ must converge to a function $\tilde{v} \in H^\RegExponent$ and to $v(-T_{lwp}^{-1})$ in $L^2$. Since $H^\RegExponent$ is dense in $L^2$, $v(-T_{lwp}^{-1})=\tilde{v}\in H^\RegExponent$. 

Since the \nls evolution and the pseudoconformal transform both preserve the $L^2$ norm, if $\|u_0\|_{L^2} < \|Q\|_{L^2}$, then $\|v(-T_{lwp}^{-1})\|_{L^2}< \|Q\|_{L^2}$. 

Therefore, in the defocusing case, from the assumption of global well-posedness in $H^\RegExponent$, $v$ extends to a function $v:\Reals\times\Reals^d\rightarrow\Complex$ with $v:\Reals\rightarrow H^\RegExponent(\Reals^d)$. For $t>0$, $u$ can be defined by $u=\Conf^{-1}[v]$. By Theorem \ref{T10.2}, this extension of $u$ is a solution to the \nls equation on $[0,\infty)$. For $t<0$ all the arguments of the paper can be reproduced to define $u$ on $(\infty,0]$. Thus, $u$ is a solution to the \nls equation, has initial data $u_0$, and is defined for all $t$. 

In the focusing case, since the \nls evolution and the pseudoconformal transform both preserve the $L^2$ norm, if $\|u_0\|_{L^2} < \|Q\|_{L^2}$, then $\|v(-T_{lwp}^{-1})\|_{L^2}< \|Q\|_{L^2}$ and the same argument can be applied with the additional $L^2$ norm hypothesis. 
\end{proof}

\begin{remark}
\label{RBackwardsInTime}
The process of taking $u$ on $[0,T_{lwp}]$, applying the pseudoconformal transform to get $v$ on $(-\infty,-T_{lwp}^{-1}]$, and then extending $v$ globally in time provides a function $v$ which is defined for positive time. We remark that there is no clear relation between $v$ at positive time and $u$ at negative time. In some sense, $v$ at positive time corresponds to the evolution of $u$ ``beyond infinite'', and, unless it is known {\it a priori} that the scattering states $u_+$ and $u_-$ satisfy
\begin{equation*}
\lim_{t\rightarrow\infty} \|e^{-\frac{i x^2}{t}}e^{it\Delta}u_+ - e^{\frac{i x^2}{t}}e^{-it\Delta}u_-\|_{L^2}=0 , 
\end{equation*}
there is no reason to believe that $v$ at positive transformed time corresponds to $u$ at negative time. 
\end{remark}

\section{A scattering lemma}
\label{s3}

\begin{lemma}
\label{ScatteringLemma}
If $u$ is a solution of the \nls equation \eqref{NLS} which exists globally and for which $v$ exists globally as well, then $\|u\|_{S^0[0,\infty)}\leq\infty$ and there exists a $u_+$ for which
\begin{equation*}
\lim_{t\rightarrow\infty} \|u(t)-e^{it\Delta}u_+\|_{L^2}=0 .
\end{equation*}
\end{lemma}
\begin{proof}
If $u$ and $v$ exist globally, then, in particular they exist for $t\in[0,1]$ and $\transT\in[-1,0]$ respectively. By the maximal time blow up Theorem \ref{MaxTimeBlowUp}, $\|u\|_{S^0([0,1])}<\infty$ and $\|u\|_{S^0([1,\infty))}=\|v\|_{S^0([-1,0])}<\infty$. Thus $\|u\|_{S^0([0,\infty))}<\infty$. By Duhamel's principle and the Strichartz estimates, this means that $\lim_{t\rightarrow\infty} \|u(t)-e^{it\Delta}u_+\|_{L^2}=0$. 
\end{proof}

\begin{remark}
\label{RBackwardScattering}
As in Remark \ref{RBackwardsInTime}, it is necessary to distinguish between the forward in time extension of the pseudoconformal transform and the backwards in time extension of the pseudoconformal transform. If $\tilde{v}$ denotes the pseudoconformal transform of $u$ at negative time, then if $u$ and $\tilde{v}$ both extend globally in time, by the same argument as in this lemma, $\|u\|_{S^0((-\infty,0])}<\infty$ and there is a $u_-$ for which $\lim_{t\rightarrow-\infty}\|u-e^{-it\Delta}u_-\|_{L^2}=0$. 
\end{remark}

\begin{corollary}
\label{ScatteringLemmaForu}
Under the hypotheses of Proposition \ref{P3.1}, there are functions $u_\pm\in \Xsigma$ such that 
\begin{equation*}
\lim_{t\rightarrow\pm\infty} \|e^{\mp it\Delta} u(t)- u_\pm\|_{\Xsigma} =0 . 
\end{equation*}
\end{corollary}
\begin{proof}
The construction in the proof of Proposition \ref{P22.1} shows that both $u$ and $v$ exist globally, and hence $u$ scatters forward in time by Lemma \ref{ScatteringLemma}. As noted in Remark \ref{RBackwardScattering}, the same occurs backwards in time. This establishes the existence of $u_\pm\in L^2$. 

We now introduce linearly advanced and retarded versions of $u$ and $v$. These have two time variables, one to record the time variable associated with the \nls evolution, and one for the advancement or retardation by the \ls evolution. 
\begin{align*}
\phi(t,t') =& e^{i(t'-t)\Delta}u(t)\\
\psi(\transT,\transT') =& e^{i(\transT'-\transT)\Delta} v(\transT)
\end{align*}
The function $\phi(t,\bullet)$ is a linear solution with initial data $u(t)$ at time $t'=t$. The function $\psi(\transT,\bullet)$ is the analogous function with initial data $v(\transT)$ at time $\transT'=\transT$. Since $v(-t^{-1})$ is the pseudoconformal transform of $u(t)$ at time $t$, and the pseudoconformal transform preserves the \ls evolution, the pseudoconformal transform of $\phi(t,\bullet)$ with respect to the spatial variable and the second time variable is $\psi(-t^{-1},\bullet)$. 

Denoting a solution to the \ls equation by $\phi$ and its pseudoconformal transform by $\psi$, it is known that $\psi(0)$ is the Fourier transform of $\phi(0)$. If $u_0\in\Xsigma$, by the construction in Section \ref{s2}, then $v(\transT)\in H^\RegExponent$. Since the \ls evolution  preserves the $H^\RegExponent$, for all $\transT'$, $\psi(\transT,\transT')\in H^\RegExponent$, and, in particular $\psi(\transT,0)\in H^\RegExponent$. This proves that the linearly retarded version of $u$ evolves in $\Xsigma$. 
\begin{align*}
e^{-it\Delta}u(t)
= \phi(t,0)
= \Fourier{\psi(\transT,0)} \in \Xsigma
\end{align*}

Since $e^{-it\Delta}u(t)$ evolves in $\Xsigma$, it is now possible to consider scattering in $\Xsigma$, despite the fact that $u(t)$ itself does not evolve in $\Xsigma$. From the $L^2$ scattering, it follows that $u_+ = \Fourier{\psi(0,0)}=\Fourier{v(0)}$. It now remains to show this limit holds in $\Xsigma$. 
\begin{align*}
\| e^{-it\Delta}u(t)-u_+ \|_{\Xsigma}
=& \| \phi(t,0) - \Fourier{v(0)} \|_{\Xsigma}\\
=& \| \psi(\transT,0) - v(0) \|_{H^\RegExponent}\\
=& \| e^{-i\transT\Delta}v(\transT) - v(0) \|_{H^\RegExponent}\\
=& \| v(\transT) - e^{i\transT\Delta}v(0) \|_{H^\RegExponent}\\
\lim_{t\rightarrow\infty} \| e^{-it\Delta}u(t)-u_+ \|_{\Xsigma}
=& \lim_{\transT \rightarrow 0^+} \| v(\transT) - e^{i\transT\Delta}v(0) \|_{H^\RegExponent} \rightarrow 0
\end{align*}

The same argument holds as $t\rightarrow -\infty$, with the usual remark on the difference between the pseudoconformal transforms for positive and negative times. 

\end{proof}

\newpage

\newcommand{\ApproxIndex}{\rho}

\newcommand{\di}{\mathrm{d}}

\begin{center}
{\sc\large\bf Errata}
\end{center}

In this paper, we consider the initial value problem for the $L^2$-critical nonlinear Schr\"odinger equation for $u:\Reals\times\Reals^d\rightarrow\Complex$: 
\begin{align}
\begin{cases}
i\partial_t u +\Delta u= \lambda|u|^\frac4d u \\
u(t_0,x)=u_0(x) ,
\end{cases}
\label{eq:NLS}
\end{align}
with initial data $u_0$ in $H^{0,s}$, the weighted $L^2$ space with norm
\begin{align*}
\|\psi\|_{H^{0,s}}&= \left( \int_{\Reals^d} (1+x^2)^s |\psi(x)|^2 \di x \right)^{1/2} .
\end{align*}
Our main result was that if equation \eqref{eq:NLS} is globally well-posed for initial data in the $H^s$ Sobolev space, then all initial data $u_0\in H^{0,s}$ generates global solutions $u\in L^\infty_t(L^2_x)$. We also proved that such solutions scatter. 

The proof presented in Section 2 approximates $u_0\in H^{0,s}$ by a sequence of regularised approximators in $H^{0,s}\cap H^1$. It has been brought to our attention that there appears to be an error in the final line of the proof of Lemma 2.1. In particular, although the growth of the $H^{0,s}$ norms of the approximators is controlled by the $H^1$ norm, it seems that the growth of the divergence of two approximators in $H^{0,s}$, $\|u'-u''\|_{H^{0,s}}$ is not controlled by the $H^{1}$ norms of the approximators. Fortunately, the same argument works if we replace $H^{0,s}\cap H^1$ by $H^{0,s}\cap H^{\ApproxIndex}$ with
$\ApproxIndex=(d+1)/2$. 

In the proof of Lemma 2.1, the application of Theorem 1.4 applies equally well to the persistence of regularity in $H^{\ApproxIndex}$ as it did in $H^1$ in the original proof. Thus, the same norm bound applies. The calculation of the growth of the $H^{0,s}$ norm of $u'$ is correct, and, since $1\leq(d+1)/2$, we can trivially replace the $H^{1}$ norm by the $H^\ApproxIndex$ norm. In the final calculation of Lemma 2.1, where there seems to be an error, we can sill use $\langle\cdot,\cdot\rangle$ for the $L^2$ inner product and $\langle x\rangle$ for $(1+x^2)^{1/2}$ and make the estimate:
\begin{align}
\frac{\di}{\di t} \|u'(t)-u''(t)\|_{H^{0,s}}^2
&=
i\langle \Delta (u'-u''),\langle x\rangle^{2s} (u'-u'')\rangle
\label{eq:RawTDeriv}\\
&\quad-i\langle (u'-u''),\langle x\rangle^{2s} \Delta (u'-u'')\rangle \nonumber\\
&\quad-i\lambda\langle |u'|^{\frac4d}u'-|u''|^\frac4d u'',\langle x\rangle^{2s} (u'-u'')\rangle\nonumber\\
&\quad+i\lambda\langle u'-u'',\langle x\rangle^{2s} (|u'|^{\frac4d}u'-|u''|^\frac4d u'')\rangle .\nonumber
\end{align}
The first two terms on the right can be estimated as in the estimate of the norm of $u'$, so that they are dominated by
\begin{align}
C\|u'-u''\|_{H^{0,s}}\|u'-u''\|_{H^1} .
\label{eq:TDerivLinearTerms}
\end{align}
The argument which needs correcting concerns the estimate on the remaining terms. These are controlled by
\begin{align*}
C|\langle |u'|^{\frac4d}u'&-|u''|^\frac4d u'',\langle x\rangle^{2s} (u'-u'')\rangle| \\
\leq&C \|\langle x\rangle^s(u'-u'')\|_{L^2} \| \langle x\rangle^s (|u'|^{\frac4d}u'-|u''|^\frac4d u'')\|_{L^2} \\
\leq&C \|u'-u''\|_{H^{0,s}}\left( \|\langle x\rangle^s u'\|_{L^2} \| u'\|_{L^\infty}^\frac4d +\|\langle x\rangle^s u''\|_{L^2} \|u''\|_{L^\infty}^\frac4d\right) \\
\leq&C \|u'-u''\|_{H^{0,s}}\left( \|u'\|_{H^{0,s}} \| u'\|_{H^\ApproxIndex}^\frac4d +\|u''\|_{H^{0,s}} \|u''\|_{H^\ApproxIndex}^\frac4d\right) .
\end{align*}
Here we have used that $\ApproxIndex>d/2$, to control the $L^\infty$ norm. By applying the earlier estimate on the growth of the $H^\ApproxIndex$ and $H^{0,s}$ norms and restricting $t\leq1$, we have the stronger bound by 
\begin{align}
C \|u'-u''\|_{H^{0,s}}
\left( (\|u'_0\|_{H^{0,s}}+\| u'_0\|_{H^\ApproxIndex}) \| u'_0\|_{H^\ApproxIndex}^\frac4d +(\|u''\|_{H^{0,s}}+\|u''_0\|_{H^\ApproxIndex}) \|u''_0\|_{H^\ApproxIndex}^\frac4d\right) .
\label{eq:TDerivNonlinearTerms}
\end{align}
Combining the estimates in \eqref{eq:TDerivLinearTerms} and \eqref{eq:TDerivNonlinearTerms}, we have that
\begin{align*}
&\frac{\di}{\di t} \|u'(t)-u''(t)\|\\
&\leq C\left( (\|u'_0\|_{H^{0,s}}+\| u'_0\|_{H^\ApproxIndex}) \left(1+\| u'_0\|_{H^\ApproxIndex}^\frac4d\right) +(\|u''\|_{H^{0,s}}+\|u''_0\|_{H^\ApproxIndex}) \left(\|u''_0\|_{H^\ApproxIndex}^\frac4d+1\right) \right) .
\end{align*}
Thus, the statement of Lemma 2.1 remains valid if we (i) replace the space $H^1$ by the space $H^\ApproxIndex$, (ii) replace the $H^1$ norm by the $H^\ApproxIndex$ norm in estimates of the growth of Sobolev norms, and (iii) replace the $\|u'\|_{H^1}$ norm by
\begin{align}
(\|u'_0\|_{H^{0,s}}+\| u'_0\|_{H^\ApproxIndex}) \left(1+\| u'_0\|_{H^\ApproxIndex}^\frac4d\right)
\label{eq:UglyNormlike}
\end{align}
and similarly for $u''$ in estimates of the growth of the $H^{0,s}$ norm. 

In Sections 2.2 and 2.3, it is sufficient to replace $H^1$ by $H^\ApproxIndex$. In the proof of 2.3 and 2.4, it is sufficient to use $H^1$ throughout and then conclude with the embedding of $H^\ApproxIndex\hookrightarrow H^1$. In Lemma 2.5, it is sufficient to replace $H^1$ by $H^\ApproxIndex$, since both are sufficient to put $v'$ and $v''$ in $H^s$. 

In Section 2.4, it is sufficient to make the same three changes as in Section 2.1. The replacement of the $H^1$ norm by the term in \eqref{eq:UglyNormlike} does not prevent the technique of taking $t\rightarrow0$ to leave dependence only on the $H^{0,s}$ norm. No further use of the regularised approximators is made in the paper. 

\vfill
\noindent{\bf Acknowledgements.} We would like to thank Se\c{c}kin
Demirba\c{s} of Bogazici University
for bringing the error in the original paper to our attention. 
\end{document}

%% file: NLSConformalPicture1.pstex_t
\begin{picture}(0,0)%
\includegraphics{NLSConformalPicture1.pstex}%
\end{picture}%
\setlength{\unitlength}{4144sp}%
\begingroup\makeatletter\ifx\SetFigFont\undefined%
\gdef\SetFigFont#1#2#3#4#5{%
  \reset@font\fontsize{#1}{#2pt}%
  \fontfamily{#3}\fontseries{#4}\fontshape{#5}%
  \selectfont}%
\fi\endgroup%
\begin{picture}(7632,4545)(76,-4684)
\put(2296,-1996){\makebox(0,0)[lb]{\smash{{\SetFigFont{12}{14.4}{\rmdefault}{\mddefault}{\updefault}{\color[rgb]{0,0,0}$u$}%
}}}}
\put(2026,-1366){\makebox(0,0)[lb]{\smash{{\SetFigFont{12}{14.4}{\rmdefault}{\mddefault}{\updefault}{\color[rgb]{0,0,0}$u''$}%
}}}}
\put( 91,-4606){\makebox(0,0)[lb]{\smash{{\SetFigFont{12}{14.4}{\rmdefault}{\mddefault}{\updefault}{\color[rgb]{0,0,0}$\Xsigma\cap H^1$}%
}}}}
\put(7021,-4111){\makebox(0,0)[lb]{\smash{{\SetFigFont{12}{14.4}{\rmdefault}{\mddefault}{\updefault}{\color[rgb]{0,0,0}$L^2$}%
}}}}
\put(4501,-4516){\makebox(0,0)[lb]{\smash{{\SetFigFont{12}{14.4}{\rmdefault}{\mddefault}{\updefault}{\color[rgb]{0,0,0}$H^\sigma$}%
}}}}
\put(6841,-2221){\makebox(0,0)[lb]{\smash{{\SetFigFont{12}{14.4}{\rmdefault}{\mddefault}{\updefault}{\color[rgb]{0,0,0}$v$}%
}}}}
\put(6076,-1726){\makebox(0,0)[lb]{\smash{{\SetFigFont{12}{14.4}{\rmdefault}{\mddefault}{\updefault}{\color[rgb]{0,0,0}$v'$}%
}}}}
\put(6616,-1726){\makebox(0,0)[lb]{\smash{{\SetFigFont{12}{14.4}{\rmdefault}{\mddefault}{\updefault}{\color[rgb]{0,0,0}$v''$}%
}}}}
\put(1486,-1411){\makebox(0,0)[lb]{\smash{{\SetFigFont{12}{14.4}{\rmdefault}{\mddefault}{\updefault}{\color[rgb]{0,0,0}$u'$}%
}}}}
\put(721,-961){\makebox(0,0)[lb]{\smash{{\SetFigFont{12}{14.4}{\rmdefault}{\mddefault}{\updefault}{\color[rgb]{0,0,0}$T_{lwp}$}%
}}}}
\put(1036,-3031){\makebox(0,0)[lb]{\smash{{\SetFigFont{12}{14.4}{\rmdefault}{\mddefault}{\updefault}{\color[rgb]{0,0,0}$t$}%
}}}}
\put(946,-3661){\makebox(0,0)[lb]{\smash{{\SetFigFont{12}{14.4}{\rmdefault}{\mddefault}{\updefault}{\color[rgb]{0,0,0}$0$}%
}}}}
\put(5401,-3526){\makebox(0,0)[lb]{\smash{{\SetFigFont{12}{14.4}{\rmdefault}{\mddefault}{\updefault}{\color[rgb]{0,0,0}$-t^{-1}$}%
}}}}
\put(5221,-1366){\makebox(0,0)[lb]{\smash{{\SetFigFont{12}{14.4}{\rmdefault}{\mddefault}{\updefault}{\color[rgb]{0,0,0}$-T_{lwp}^{-1}$}%
}}}}
\put(1936,-4201){\makebox(0,0)[lb]{\smash{{\SetFigFont{12}{14.4}{\rmdefault}{\mddefault}{\updefault}{\color[rgb]{0,0,0}$L^2$}%
}}}}
\put(2206,-3526){\makebox(0,0)[lb]{\smash{{\SetFigFont{12}{14.4}{\rmdefault}{\mddefault}{\updefault}{\color[rgb]{0,0,0}$\delta$}%
}}}}
\put(4006,-2761){\makebox(0,0)[lb]{\smash{{\SetFigFont{12}{14.4}{\rmdefault}{\mddefault}{\updefault}{\color[rgb]{0,0,0}{\Large $\Conf$}}%
}}}}
\end{picture}%